\documentclass[oneside]{amsart}
\usepackage{epsfig,amsmath,latexsym,amssymb,amsfonts,amsthm,cite}
\usepackage{graphicx}
\usepackage{graphics}
\usepackage{color,bm}
\graphicspath{{../img/}}
\usepackage[all]{xy}
\newtheorem{theorem}{Theorem}[section]
\newtheorem{lemma}[theorem]{Lemma}
\newtheorem{proposition}[theorem]{Proposition}

\newtheorem{question}{Question}
\newtheorem{definition}{Definition}

\newtheorem{corollary}[theorem]{Corollary}
\newtheorem{main}{Theorem}

\def\F{\mathcal{F} }

\def\Q{\mathbb{Q} } 
\def\R{\mathbb{R} } 
\def\Z{\mathbb{Z} } 
\def\nbd{neighborhood } 
\def\nbds{neighborhoods } 
\def\R{\mathbb{R} }

\def\-{\ominus} 
\def\+{\oplus} 
\def\0{\circ}

\usepackage{color,bm}

\author{Tomoo Yokoyama}
\date{\today}
\address{Department of Mathematics, Saitama University, Saitama, 255 Shimookubo, Sakura Ward, Saitama, 338-8570, Japan\\}
\email{tyokoyama@rimath.saitama-u.ac.jp}
\thanks{The author was partially supported by JSPS Grant Number 21H00980}
\makeatletter
\@namedef{subjclassname@2020}{%
  \textup{2020} Mathematics Subject Classification}
\makeatother

\title[Recurrence for semi-decompositions]{Recurrence for semi-decompositions}



\begin{document}


\maketitle

\begin{abstract}
This paper constructs a foundation to analyze semi-group actions, group actions, filtrations, and decompositions in a unified manner. In fact, though the studies of decomposition can be applied to foliated spaces and group actions, they can not be applied to semi-group actions and filtrations in general because filtration and the set of orbits of a semi-group need not be decompositions of the base spaces. To analyze these concepts in a unified manner, we introduce a concept of a semi-decomposition which is a natural generalization of these concepts because similar relations among recurrence and their variants to group actions and decompositions hold for semi-decompositions. On the other hand, we demonstrate the difference between the recurrent concepts for group actions and those even for semi-group actions.
\end{abstract}

\maketitle

\section{Introduction}

Recurrence is one of the most fundamental properties of dynamics. It is generalized into several concepts (e.g. equicontinuity, variants of almost periodicity, distal property, characteristic $0$ property) for actions of groups and semi-groups. Some concepts for group actions also are introduced for decomposition (e.g. foliated spaces). Though the studies of decomposition can be applied to foliated spaces and group actions, they can not be applied to semi-group actions in general because filtration and the set of orbits of a semi-group need not be decompositions of the base spaces. This paper introduces semi-decomposition, which is a generalization of the concepts, to analyze semi-group actions, group actions, filtrations, and decompositions in a unified manner. 

From a semi-group action point of view, Gottschalk characterizes almost periodicity \cite{Gottschalk1946almost}, Ellis, Ellis, and Nerurkar study the recurrence in this century \cite{Ellis2001top}, and Auslander and Dai introduce variants of recurrence for semi-group actions and generalize several results to them for flows to semi-group actions in recent decades \cite{Auslander2019mini}. 
In this paper, we develop various relations among variants of recurrence into those for semi-decompositions. 

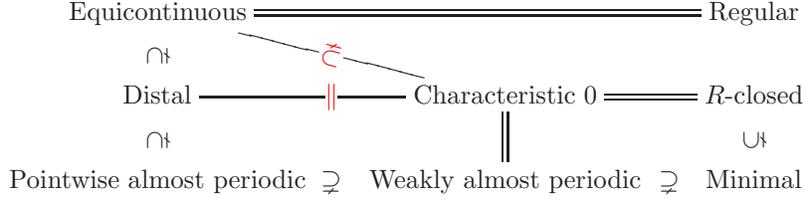
\begin{figure}[t]
\[
\xymatrix@=18pt{
\text{Equicontinuous} \ar@{-}[dr]|{\rotatebox{165}{{\color{red}$\supsetneq$}}} \ar@{}[d]|{\rotatebox{90}{$\supsetneq$}} \ar@{=}[rr] &   & \text{Regular}  \\
\text{Distal}  \ar@{}[d]|{\rotatebox{90}{$\supsetneq$}} \ar@{-}[r]|{\rotatebox{0}{{\color{red}$\parallel$}}}  &  \text{Characteristic } 0  \ar@{=}[d]  & R\text{-closed} \ar@{=}[l] 
 \\
 \text{Pointwise almost periodic} & \text{Weakly almost periodic} \ar@{}[l]|{\rotatebox{0}{$\supsetneq$}} & \text{Minimal} \ar@{}[l]|{\hspace{30pt}\rotatebox{0}{$\supsetneq$}} \ar@{}[u]|{\rotatebox{90}{$\subsetneq$}} 
  }
\]
\caption{Relations among concepts for group actions on compact Hausdorff spaces.}
\label{Fig:implications}
\end{figure}

From a group action point of view, there are various relations among variants of recurrence (see Figure~\ref{Fig:implications}). 
For instance, the equivalence of equicontinuity and closedness of the orbit closure relation (i.e. $R$-closedness) for the product action for group actions on compact Hausdorff spaces is shown \cite[Lemma~4.3 p.25 and Proposition~4.11 p.29]{ellis1969lectures}. 
On the other hand, the equivalence of equicontinuity and the closedness of the orbit closure relation does not hold \cite{knight1986cer,salem2018group}. 
The equivalence of closedness of the orbit closure relations and characteristic $0$ property is shown for group actions on Hausdorff spaces in \cite[Theorem~9]{knight1986cer} (also \cite[Theorem~1.7]{elaydi1982char}) and more generally for decompositions on Hausdorff space \cite[Corollary 3.1]{yokoyama2013recurrence}. 
Any equicontinuous group action on a Hausdorff space is distal (in the sense of \cite[Definition~4.5 p.41]{ellis2014automorphisms}) (cf. \cite[Proposition~15.3]{ellis2014automorphisms}).
The distal property implies pointwise almost periodicity (i.e. $T_1$-separation axiom for the orbit class space) for any group action on a locally compact Hausdorff space such that any orbit closures are compact, from \cite[Theorem~1]{ellis1958distal} and \cite[Theorem~7 p.11]{auslander1988minimal}. 
%
%
%
This paper demonstrates similar relations for semi-decompositions and the difference between the recurrent concepts for group actions and those even for semi-group actions. 

In addition, the equicontinuty implies characteristic $0$ properties and the closedness of the orbit closure relation for a group action on a compact Hausdorff space \cite[Proposition~4.10]{ellis1969lectures}, for abelian group action on locally compact Hausdorff uniformizable space \cite[Theorem~4.2]{elaydi1982p}, and group actions with compact orbit closure on locally compact Hausdorff uniformizable spaces \cite[Theorem~2, Proposition~6, Theorem~9]{knight1986cer}. 
Moreover, the implication is stated in \cite[Exercise~6(a) p.46]{auslander1988minimal} 
in the group action case.  
%
Therefore we ask whether this implication holds in general. 
More precisely, we ask the following question. 

\begin{question}\label{q:01}
Does the equicontinuity make the orbit closure relations closed in general?
\end{question}

In this paper, we give an affirmative answer for Question~\ref{q:01} in the case of semi-group actions on metrizable spaces (Corollary~\ref{lem:equi001}), and a similar statement holds for foliated spaces (see Corollary~\ref{th:fol}). 
More generally, we have the following relation. 


\begin{main}\label{th:equi_decomp}
Every weakly equicontinuous semi-decomposition on a metrizable space is $R$-closed and of characteristic $0$. 
\end{main}

%
Equicontinuity does not imply pointwise almost periodicity for semi-group actions (see examples in \S\ref{sec:ex_analytic}--\ref{sec:ex_r_nonec} and \S\ref{ex:r_non_pap}). 
On the other hand, we give the following implication in the case of decomposition on metrizable spaces. 

\begin{main}\label{th:equi}
The following statements hold for any weakly equicontinuous decomposition on a metrizable space: 
\\
{\rm(1)} The decomposition is of characteristic $0$. 
\\
{\rm(2)} The decomposition is $R$-closed.  
\\
{\rm(3)} The class space is Hausdorff.  
\\
Moreover, if any element closures are compact, then the following statements hold: 
\\
{\rm(4)} The decomposition is weakly almost periodic. 
\\
{\rm(5)} The class decomposition is an upper semi-continuous decomposition.  
\\
{\rm(6)} For any neighborhood $U$ of an element closure, there is an open invariant neighborhood $V \subseteq U$ of the element closure. 
\end{main}

Notice that the previous theorem implies the result for group actions (Corollary~\ref{lem:equi2R_metric}). 
%
%
As mentioned above, every equicontinuous group action on a compact Hausdorff space implies distal property.
Therefore we also ask the following question, which Hiraku Nozawa poses. 

\begin{question}\label{q:02}
Is there any relation between $R$-closedness {\rm(}the closedness of the orbit closure relation{\rm)} and the distal property for group actions?
\end{question}

\begin{figure}[t]
\[
\xymatrix@=18pt{
\text{Equicontinuous}  \ar@{}[r]|{\rotatebox{0}{{\color{red}$\subsetneq$}}} &\text{Weakly equicontinuous} \ar@{}[d]|{\rotatebox{90}{{\color{red}$\supsetneq$}}} \ar@{}[dr]|{\rotatebox{160}{{\color{red}$\supsetneq$}}} &  \\
\text{Distal}  &  \text{Characteristic } 0 \ar@{}[d]|{\rotatebox{90}{{\color{red}$\subsetneq$}}} & R\text{-closed} \ar@{}[l]|{\rotatebox{0}{{\color{red}=}}}  \\
 \text{Pointwise almost periodic} & \text{Weakly almost periodic}  \ar@{}[l]|{\rotatebox{0}{$\supsetneq$}} & \text{Minimal} \ar@{}[l]|{\hspace{30pt}\rotatebox{0}{$\supsetneq$}} \ar@{}[u]|{\rotatebox{90}{$\subsetneq$}} 
  }
\]
\caption{Relations among concepts for semi-groups on metrizable spaces.}
\label{Fig:implications02}
\end{figure}
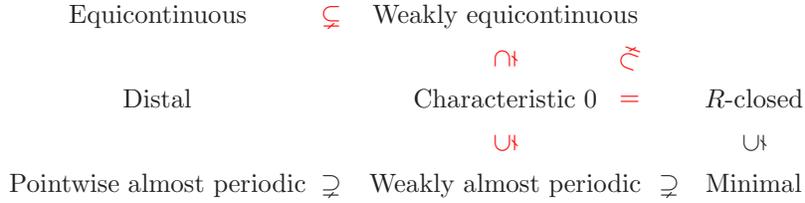

We give a negative answer for Question~\ref{q:02} as Figure~\ref{Fig:implications} even group actions on compact manifolds by constructing examples. 
Moreover, some incomparable relations hold for semi-group actions as in Figure~\ref{Fig:implications02} and more generally semi-decompositions as in Figure~\ref{Fig:implications03}. 
In addition, we have the following incomparability for semi-group actions. 

\begin{main}\label{th:incomparable_semi}
The following statements hold for semi-group actions on compact manifolds: 
\\
{\rm(1)} Distal property and $R$-closedness are incomparable even for group actions. 
\\
{\rm(2)} Distal property and weak equicontinuty are incomparable. 
\\
{\rm(3)} Pointwise almost periodicity is incomparable to each of distal property, $R$-closedness, weak equicontinuity, and equicontinuity. 
\\
{\rm(4)} Weak almost periodicity is incomparable to each of distal property, weak equicontinuity, and equicontinuity. 
\\
{\rm(5)} Minimality is incomparable to each of distal property, weak equicontinuity, and equicontinuity. 
\end{main}

Note that all incomparable relations except one between equicontinuity and distal property in Figure~\ref{Fig:implications02} are shown in the previous theorem. 
Moreover, Figures~\ref{Fig:implications}--\ref{Fig:implications03} demonstrate the difference between the concepts for group actions and those for semigroup actions.


%

The present paper consists of four sections.
In the next section, as preliminaries, we recall fundamental concepts and introduce semi-decomposition for a unified discussion for decompositions.
In \S~\ref{sec:main}, we define variants of recurrence for semi-decomposition. 
Theorem~\ref{th:equi_decomp} and Theorem~\ref{th:equi} are demonstrated. 
In the final section, we demonstrate Theorem~\ref{th:incomparable_semi} by constructing examples.

\section{Preliminaries}

\subsection{Fundamental concepts} 

We recall several fundamental concepts. 

\subsubsection{Decomposition}
By a {\bf decomposition}, we mean a family $\mathcal{F}$ of pairwise disjoint nonempty subsets of a set $X$ such that $X = \bigsqcup \mathcal{F}$, where $\bigsqcup$ denotes a disjoint union. 
Denote by $\mathcal{F}(x)$ the element of $\mathcal{F}$ containing $x$. 
For a subset $A \subseteq X$, set $\F(A) : = \bigcup_{x \in A} \F(x)$, called the {\bf saturation} of $A$. 
The closure of the element containing a point is called the {\bf element closure} of the point. 
The {\bf decomposition space} $\bm{X/\F}$ of $X$ is a quotient space $X/\sim_T$ defined by $x \sim_{\F} y$ if $\F(x) = \F(y)$.
Since connectivity is not required, the set of orbits of group actions are also decompositions.

\subsubsection{Separation axioms and relative concepts}

A point $x$ of a topological space $X$ is $\bm{T_0}$ (or {\bf Kolmogorov}) if for any point $y \neq x \in X$ there is an open subset $U$ of $X$ such that $|\{x, y \} \cap U| =1$, where $|A|$ is the cardinality of a subset $A$.
A topological space is $\bm{T_0}$ (or {\bf Kolmogorov}) if each point is $T_0$.
For a topological space $X$, define the {\bf class} $\bm{\hat{x}} := \{ y \in X \mid \overline{\{ x \}} = \overline{\{ y \}} \}$ for any point $x$ and a decomposition $\hat{X} := \{ \hat{z} \mid z \in X \}$ of classes. 
Here $\overline{A}$ is the closure of a subset $A$. 
Then the decomposition $\hat{X}$ is a $T_0$ space as a quotient space, which is called the {\bf $\bm{T_0}$-tification} (or Kolmogorov quotient) of $X$.

A topological space $X$ is $\bm{T_1}$ if any singletons are closed. 
Here a {\bf singleton} is a set with exactly one element. 
A topological space $X$ is regular if for any closed subset $F \subseteq X$ and any point $x \in X - F$ there are a disjoint neighborhood of $F$ and a neighborhood of $x$, and a topological space is $T_3$ if it is $T_1$ and regular. 
Here $B - C$ is used instead of the set difference $B \setminus C$ when $B \subseteq C$.
A topological space $X$ is completely regular if for any closed subset $A \subseteq X$ and any point $ x \in X - A$ there is a continuous function $f \colon X \to \mathbb{R}$ with $f(x)=1$ and $f\vert_{A}=0$. 
A topological space $X$ is {\bf Tychonoff} (or $T_{3\frac{1}{2}}$) if it is Hausdorff and completely regular. 
A topological space $X$ is {\bf uniformizable} if there is a uniformity on $X$ that induces the topology of $X$. 
In other words, a topological space is uniformizable if and only if it is homeomorphic to a uniform space equipped with the topology of the uniformity. 
It is known that any topological space is uniformizable if and only if it is completely regular (cf. \cite{sen2010causality}), and so that any Hausdorff space is uniformizable if and only if it is Tychonoff. 
A topological space is normal if any pair of two disjoint closed subsets has disjoint open neighborhoods, and a topological space is $T_4$ if it is $T_1$ and normal.

\subsubsection{Hausdorff distance}

Let $(X,d)$ be a metric space. 
The {\bf Hausdorff distance} $d_{\rm {H}}: (2^X - \{ \emptyset \}) \times 2^X - \{ \emptyset \}) \to  \R_{\geq 0} \sqcup \{ \infty \}$ is an extended quasimetric defined as follows: for any non-empty subsets $A, B \subseteq X$, 
\[
d_{H}(A,B) := \inf\{\varepsilon > 0 \mid A \subseteq U_{\varepsilon }(B), B \subseteq U_{\varepsilon }(A) \}
\]
where $U_{\varepsilon }(C) := \bigcup _{x\in C} \{ y \in X \mid d(y,x) < \varepsilon \}$.  
    
\subsubsection{Pre-order}
A binary relation $\leq$ on a set $P$ is a {\bf pre-order} if it is reflexive (i.e. $a \leq a$ for any $a \in P$) and transitive (i.e. $a \leq c$ for any $a, b, c \in P$ with $a \leq b$ and $b \leq c$).
For a point $x$ of a pre-ordered set $P$, define the {\bf downward closure} $\mathop{\downarrow}x := \{ y \in P \mid y \leq x \}$ of the point $x$. 

\subsection{Concepts for semi-decompositions}

\subsubsection{Semi-decomposition}

Though the studies of decomposition can be applied to foliated spaces and group actions, they can not be applied to semi-group actions in general because the set of orbits of a semi-group need not be a decomposition of the base space. 
Therefore we introduce a concept of semi-decomposition to analyze semi-group actions and decompositions in a unified manner as follows.

\begin{definition}
A mapping $\F \colon X \to 2^X$ from a set $X$ is a {\bf semi-decomposition} if it satisfies the following two conditions: 
\\
{\rm(1)} $x \in \F(x)$ for any $x \in X$. 
\\
{\rm(2)} $\F(x) \subseteq \F(y)$ for any $x, y \in X$ with $x \in \F(y)$. 
\end{definition}
The closure of the element containing a point is called the {\bf element closure} of the point. 

By definition, any decompositions are semi-decompositions. 
Note that the set of orbits of any semi-group action is a semi-decomposition and that the set of simplices of the simplicial complex has a semi-decomposition structure. 
Moreover, any filtrations of topological spaces and any stratifications on topological spaces have semi-decomposition structures. 
In particular, the set of (connected components of) sub-level sets of any function is a semi-decomposition. 
Here a {\bf stratification} on a topological space $X$ is a continuous mapping $\pi \colon X \to P$ from $X$ to a poset $P$ with the Alexandroff topology (cf. \cite[Definition A.5.2]{luriehigher} and \cite[Definition 2.1.3]{Ayala2017local}), and the subsets $\F(x) := \{ \bigcup_{y \leq x} \pi^{-1}(y) \mid x \in P \}$ for any points $x \in X$ form a semi-decomposition of the stratification. 
Notice that a semi-decomposition on a topological space has no stratification in general, even for decompositions generated by groups. 
In fact, the set of orbits of an irrational rotation on a circle, which is a minimal diffeomorphism, is a decomposition but does not have canonical stratifications because any orbits are not closed but maximal elements in the decomposition with respect to the inclusion order. 

From now on, we use $\F$ for a semi-decomposition on a topological space $X$ unless otherwise stated. 

\begin{figure}[t]
\[
\xymatrix@=18pt{
 &\text{Weakly equicontinuous} \ar@{}[d]|{\rotatebox{90}{{\color{red}$\supsetneq$}}} \ar@{}[dr]|{\rotatebox{160}{{\color{red}$\supsetneq$}}}  & \\
 &  \text{Characteristic } 0 \ar@{}[d]|{\rotatebox{90}{{\color{red}$\subsetneq$}}} & R\text{-closed} \ar@{}[l]|{\rotatebox{0}{{\color{red}=}}}  \\
 \text{Pointwise almost periodic} & \text{Weakly almost periodic}  \ar@{}[l]|{\rotatebox{0}{$\supsetneq$}} & \text{Minimal} \ar@{}[l]|{\hspace{30pt}\rotatebox{0}{$\supsetneq$}} \ar@{}[u]|{\rotatebox{90}{$\subsetneq$}} 
  }
\]
\caption{Relations among concepts for semi-decomposition on metrizable spaces.}
\label{Fig:implications03}
\end{figure}
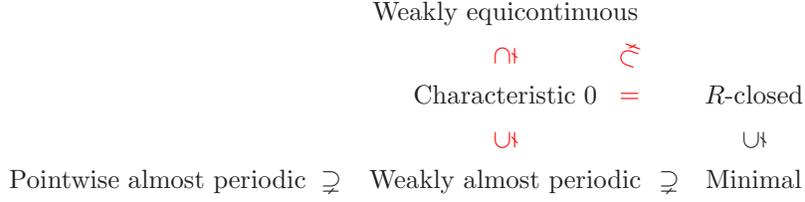

In this paper, we show the relations of properties for semi-decompositions as in Figure~\ref{Fig:implications03} (see Lemmas~\ref{lem:equiv_ch_r}--\ref{th:mini}
and an example in \S\ref{sec:ec_npap}). 

\subsubsection{A one-to-one correspondence between semi-decompositions and pre-orders}

There is a one-to-one correspondence between semi-decompositions and pre-orders as follows: 

For a semi-decomposition $\F$ on a set $X$, define a pre-order $\leq$ on $X$ by $x \leq y$ if $x \in \F(y)$ (equivalently $\F(x) \subseteq \F(y)$). 
In other words, the induced pre-order is the order induced by the inclusion order on the semi-decomposition $\F$. 

For a pre-order $\leq$ on a set $X$, define a semi-decomposition $\F$ by $\F(y) := \mathop{\downarrow} y = \{ x \in X \mid x \leq y \}$ (i.e. $\F = \{ \mathop{\downarrow} y \mid y \in X \}$). 
In other words, the induced semi-decomposition is the set of downward closures of points.

\subsubsection{Properties of semi-decomposition}
For any point $x \in X$, define the {\bf element class} $\hat{\mathcal{F}}(x)$ as follows:
$$\hat{\mathcal{F}}(x) := \{ y \in X \mid \overline{\mathcal{F}(x)} = \overline{\mathcal{F}(y)} \} \subseteq \overline{\mathcal{F}(x)}$$
We call $\hat{\mathcal{F}}(x)$ the element class of $x$. 
If the class semi-decomposition is a decomposition, then the decomposition is denoted by $\hat{\F}$ and is called the class decomposition. 
For the class decomposition $\hat{\F}$, the {\bf (element) class space} $X/\hat{\mathcal{F}}$ of $X$ is a quotient space $X/\sim_{\hat{\mathcal{F}}}$ defined by $x \sim_{\hat{\mathcal{F}}} y$ if $\hat{\mathcal{F}}(x) = \hat{\mathcal{F}}(y)$ (or equivalently $\overline{\mathcal{F}(x)} = \overline{\mathcal{F}(y)}$).
Note that the decomposition class space $X/\hat{\mathcal{F}}$ of a decomposition is the $T_0$-tification of the decomposition space $X/\F$.

For a semi-decomposition $\F$ on a topological space $X$, the {\bf (bilateral) prolongation} $D(x)$ of a point $x \in X$ is defined as follows: 
\[
D(x) = \{ y \in X \mid y_{\alpha} \in  \F(x_{\alpha}), y_{\alpha} \to y,  \text{ and } x_{\alpha} \to x \text{ for some nets } (y_{\alpha}),  (x_{\alpha}) \subseteq X \}
\]

We have the following observation.
 
\begin{lemma}\label{lem:prolong}
For any point $x$ of a semi-decomposition on a topological space, the element closure of $x$ is contained in the bilateral prolongation $D(x)$. 
\end{lemma}

\begin{proof}
Let $\F$ be a semi-decomposition on a topological space $X$. 
Fix a point $x \in X$. 
Since $X$ is a metric space, by $D(x) = \{ y \in X \mid y_{\alpha} \in  \F(x_{\alpha}), y_{\alpha} \to y,  \text{ and } x_{\alpha} \to x \text{ for some nets } (y_{\alpha}),  (x_{\alpha}) \subseteq X \}$, we have $\overline{\F(x)} = \{ y \in X \mid y_{\alpha} \in  \F(x), y_{\alpha} \to y,  \text{ for some net } (y_{\alpha}) \subseteq X \} \subseteq D(x)$. 
\end{proof}

\subsubsection{Properties of semi-decompositions}
A subset of $X$ is {\bf invariant} with respect to a semi-decomposition $\F$ if it is a union of elements of $\F$. 
By definition of semi-decomposition, a subset $A$ is invariant if and only if $A = \bigcup_{x \in A} \F(x)$. 
An invariant subset $A$ of $X$ is {\bf minimal} with respect to $\F$ if $A = \overline{\F(x)}$ for any point $x \in A$. 

\section{Relations of variants of recurrence}\label{sec:main}

In this section, we study relations of variants of recurrence for decompositions and semi-decompositions. 

\subsection{Characteristic $0$ property and weak equicontinuouity}

We define pointwise almost periodicity as follows. 

\begin{definition}
A semi-decomposition $\F$ on a topological space $X$ is {\bf pointwise almost periodic} if any element closures of $\F$ are minimal. 
\end{definition}
Notice that A semi-decomposition on a topological space $X$ is pointwise almost periodic if and only if the element closure relation $R = \{ (x,y) \in X \times X \mid y \in \overline{\F(x)} \}$ is symmetric.

\begin{definition}
A semi-decomposition on a topological space $X$ is {\bf $R$-closed} {\rm(}or {\bf the element closure relation is closed}{\rm)} if the element closure relation $R = \{ (x,y) \in X \times X \mid y \in \overline{\F(x)} \}$ is closed as a subset of the product space. 
\end{definition}

By \cite[Corollary~1.4]{yokoyama2013recurrence}, the $R$-closedness for decompositions implies the pointwise almost periodicity. 
Therefore the closedness of the element closure relation for decompositions implies the symmetric property. 
However, the closedness of the element closure relation for semi-decompositions does implies the symmetric property. 
In fact, there is an $R$-closed semi-decomposition that is not pointwise almost periodic such that the element closure relation is not symmetric (see an example in \S~\ref{ex:r_non_pap}). 
Note that we use the element closure relation is the inverse relation of one in \cite{yokoyama2013recurrence}. 

\subsubsection{Characteristic $0$ property}

The characteristic $0$ property is introduced by Knight \cite{knight1972char0} based on the characteristic $0^+$ property introduced by Ahamd \cite{Ahmad1970char0}. 
Consider a semi-decomposition $\F$ on a topological space $X$. 
\begin{definition}
A point $x \in X$ is of {\bf characteristic $0$} if $\overline{\F(x)} =D(x)$. 
\end{definition}
A semi-decomposition on a topological space $X$ is of {\bf characteristic 0} if so is any point of $X$. 
The characteristic $0$ property is also called D-stable property in \cite{Athanassopoulos1997stable}.
By definition, the following are equivalent for a decomposition: (1) $y \in D(x)$; (2) $x \in D(y)$. 
However, the equivalent does not hold for semi-decomposition, even semi-group action cases (see an example in \S\ref{ex:r_non_pap}). 

As mentioned above, the characteristic $0$ property and $R$-closedness coincide for decomposition on Hausdorff space. 
Moreover, the coincidence also holds for semi-decompositions, which follows from the same argument of the proof of \cite[Theorem~1.7]{elaydi1982char} for group actions on Hausdorff spaces. 

\begin{lemma}\label{lem:equiv_ch_r}
A semi-decomposition on a Hausdorff space is of characteristic $0$ if and only if it is $R$-closed. 
\end{lemma}

For the self-containedness, we state the following proof for semi-decompositions. 

\begin{proof}
Let $\F$ be a semi-decomposition on a Hausdorff space $X$ and $R$ the element closure relation. 

Suppose that $\F$ is of characteristic $0$. 
Fix a net $(x_{\alpha}, y_{\alpha})_{\alpha \in \Lambda}$ converging to a point $(x_\infty, y_\infty)$. 
By definition of $R$ and characteristic $0$ property, we obtain $y_\alpha \in \overline{\F(x_\alpha)} = D(x_\alpha)$ for any $\alpha \in \Lambda$. 
For any $\alpha \in \Lambda$, there is nets $(x_{\beta_\alpha})_{\beta_\alpha \in \Lambda_\alpha}$ and  $(y_{\beta_\alpha})_{\beta_\alpha \in \Lambda_\alpha}$ such that $y_{\beta_\alpha} \in  \F(x_{\beta_\alpha}), y_{\beta_\alpha} \to y_\alpha,  \text{ and } x_{\beta_\alpha} \to x_\alpha$. 
Consider the product partial order on $L := \Lambda \times \prod_{\alpha \in \Lambda} \Lambda_\alpha$. 
For $A = (\alpha', (\beta_\alpha)) \in L$, set $x_A := x_{\beta_{\alpha'}}$ and $y_A := y_{\beta_{\alpha'}}$.  
Then $(x_A)_{A \in L}$ and $(y_A)_{A \in L}$ are nets with $y_A \in \F(x_A)$. 

We claim that $(x_{A})_{A \in L}$ is a net converging to $x_\infty$ and that $(y_A)_{A \in L}$ is a net converging to $y_\infty$. 
Indeed, fix an open \nbd $U$ of $x_\infty$. 
By $x_{\alpha} \to x_\infty$, there is an element $\alpha' \in \Lambda$ such that $x_\alpha \in U$ for any $\alpha \geq \alpha'$. 
Then the open subset $U$ is also an open \nbd of $x_\alpha$ for any $\alpha \geq \alpha'$. 
For any $\alpha \in \Lambda$, there is an element $\beta'_\alpha \in \Lambda_\alpha$ such that $x_{\beta_\alpha} \in U$ for any $\beta_\alpha \geq \beta'_\alpha$. 
Set $A' := (\alpha', (\beta'_\alpha)) \in L$. 
Then  $x_{A} \in U$ for any $A \geq A'$. 
This means that the net $(x_A)_{A \in L}$ is a net converging to $x_\infty$.  
From the same argument for $y_\alpha$, the net $(y_A)_{A \in L}$ is a net converging to $y_\infty$. 

By definition of bilateral prolongation and characteristic $0$ property, we obtain $y_\infty \in D(x_\infty) = \overline{\F(x_\infty)}$. 
This implies that $(x_\infty, y_\infty) \in R$.

Conversely, suppose that $\F$ is $R$-closed. 
Fix a point $x \in X$. 
By Lemma~\ref{lem:prolong}, we have $\overline{\F(x)} \subseteq D(x)$. 
Choose a point $y \in D(x)$. 
Then there are nets $(y_{\alpha})$ and $(x_{\alpha})$ with $y_{\alpha} \in  \F(x_{\alpha})$, $y_{\alpha} \to y$, and $x_{\alpha} \to x$. 
By definition of the element closure relation $R$, we have $(x_\alpha, y_\alpha) \in R$ for any $\alpha \in \Lambda$. 
The closedness of $R$ implies that $(x, y) \in R$ and so $y \in \overline{\F(x)}$. 
Therefore $\overline{\F(x)} = D(x)$. 
\end{proof}

\subsubsection{Weak equicontinuouity}
Now, we introduce weakly equicontinuity for semi-decomposition as follows. 
\begin{definition}
A semi-decomposition $\F$ on a metric space $(X,d)$ is {\bf weakly equicontinuous} if  
for any $\varepsilon \in \R_{>0}$ there is a positive number $\delta \in \R_{>0}$ such that $d_H(\F(x), \F(y)) < \varepsilon$ for any points $x, y \in X$ with $d(x,y) < \delta$. 
\end{definition}

A semi-decomposition on a metrizable space $X$ is {\bf weakly equicontinuous} if there is a compatible metric $d \colon X \times X \to \R$ such that the metric space $(X,d)$ is weakly equicontinuous. 
Because semi-decompositions do not have local geometric structures in general, this definition of weak equicontinuity for semi-decompositions is sightly different from one for pseudogroups of local homeomorphisms and foliated spaces (see \cite[Definition~8.1]{Lopez2009equiconti} for details). 
However, any Riemannian foliations, which are introduced by Reinhart \cite{Reinhart1959fol}, on paracompact manifolds are weakly equicontinuous decompositions. 
We have the following relation. 

\begin{lemma}\label{th:we_ch}
Any weakly equicontinuous semi-decomposition on a metrizable space is of characteristic $0$. 
\end{lemma}

\begin{proof}
Let $\F$ be a weakly equicontinuous semi-decomposition on a metric space $(X,d)$. 
Fix a point $x \in X$. 
By Lemma~\ref{lem:prolong}, we have $\overline{\F(x)} \subseteq D(x)$. 
We claim that $\overline{\F(x)} \supseteq D(x)$. 
Indeed, fix a point $y \in D(x)$ and a positive number $\varepsilon > 0$. 
By definition of bilateral prolongation, there are nets $(y_{\alpha})$ and $(x_{\alpha})$ with $y_{\alpha} \in  \F(x_{\alpha})$ , $y_{\alpha} \to y$, and $x_{\alpha} \to x$. 
By weakly equicontinuity, there is a number $\delta \in \R_{>0}$ such that $d_H(\F(x'), \F(y')) < \varepsilon/2$ for any points $x', y' \in X$ with $d(x',y') < \delta$. 
From $y_{\alpha} \to y$ and $x_{\alpha} \to x$, there is $\alpha$ such that $d(x, x_\alpha) < \delta$ and $d(y, y_\alpha) < \delta$. 
Then $d_H(\F(x), \F(x_\alpha)) < \varepsilon/2$ and $d_H(\F(y_\alpha), \F(y)) < \varepsilon/2$. 
Since $\F(y_\alpha) \subseteq \F(x_\alpha)$, we have that $y \in \F(y) \subseteq U_{\varepsilon/2}(\F(y_\alpha)) \subseteq U_{\varepsilon/2}(\F(x_\alpha))$. 
On the other hand, we obtain $\F(x_\alpha) \subseteq U_{\varepsilon/2}(\F(x))$ and so $U_{\varepsilon/2}(\F(x_\alpha)) \subseteq U_{\varepsilon}(\F(x))$. 
Therefore $y \in U_{\varepsilon/2}(\F(x_\alpha)) \subset U_{\varepsilon}(\F(x))$ and so $y \in \overline{\F(x)}$. 
This means that $\overline{\F(x)} \supseteq D(x)$. 
Thus $\overline{\F(x)} = D(x)$. 
\end{proof}

Lemma~\ref{lem:equiv_ch_r} and Lemma~\ref{th:we_ch} imply Theorem~\ref{th:equi_decomp}. 
In particular, the following statement holds.

\begin{corollary}\label{lem:equi001}
Every equicontinuous semi-group action on a metrizable space is $R$-closed and of characteristic $0$. 
\end{corollary}

%
%


\begin{corollary}\label{th:fol}
Every weakly equicontinuous foliated space is $R$-closed and of characteristic $0$. 
\end{corollary}


\subsection{Upper-semicontinuity and weak almost periodicity}

We recall (weak) upper semi-continuity. 
\begin{definition}
A semi-decomposition $\mathcal{F}$ on a topological space $X$ is {\bf weakly upper semi-continuous} {\rm(cf.} \cite[Definition in p.9]{daverman1986decompositions}{\rm)} if each element of $\mathcal{F}$ is closed, and for any $L \in \mathcal{F}$ and any open neighborhood $U$ of $L$ there is an open invariant neighborhood $V$ of $L$ contained in $U$. 
\end{definition}

A weakly upper semi-continuous semi-decomposition $\mathcal{F}$ on a topological space $X$ is {\bf upper semi-continuous} (cf. \cite[Definition in p.13]{daverman1986decompositions}) if each element of $\mathcal{F}$ is compact. 
Notice that $R$-closedness for decompositions on compact Hausdorff spaces and the upper semi-continuity of the class decomposition coincide \cite[Theorem A]{yokoyama2017low}. 

\subsubsection{Weak almost periodicity in the sense of Gottschalk and Mason}

We recall weak almost periodicity in the sense of Gottschalk and Mason. 
\begin{definition}
A pointwise almost periodic semi-decomposition $\F$ on a topological space $X$ is {\bf weakly almost periodic} in the sense of Gottschalk and Mason \cite{mason1973weakly} if the union $\bigcup_{x\in A} \overline{\F(x)}$ of element closures for any closed subset $A$ of $X$ is closed. 
\end{definition}
In other words, a pointwise almost periodic semi-decomposition $\F$ on a topological space $X$ is weakly almost periodic  if and only if the canonical quotient map $q: X \to X/\hat{\F}$ is closed. 
Notice that weak equicontinuity does not imply weak almost periodicity for semi-group actions (see an example in \S\ref{sec:ec_npap}). 
On the other hand, we have the following statement, which is stated in \cite[Proposition 1.1.1, p.8]{daverman1986decompositions} 
when the base space is Hausdorff. 

\begin{lemma}\label{lem:wap_r03}
The following are equivalent for a semi-decomposition on a topological space: 
\\
{\rm(1)} The semi-decomposition is weakly almost periodic. 
\\
{\rm(2)} The semi-decomposition is pointwise almost periodic and the class decomposition is weakly upper semi-continuous. 
\\
In any case, every element class corresponds to the element closure. 
\end{lemma}

\begin{proof}
Let $\F$ be a semi-decomposition on a topological space $X$. 
Suppose that $\F$ is weakly almost periodic. 
By definition, the semi-decomposition $\F$ is pointwise almost periodic and so the class semi-decomposition $\hat{\F}$ is a decomposition with $\hat{\F} = \{\overline{\F(x)} \mid x \in X\}$. 
Fix a point $x \in X$ and an open \nbd $U$ of $\hat{\F}(x)$. 
Then the complement $F := X - U$ is closed. 
Since $\F$ is weakly almost periodic, the saturation $\hat{\F}(F)$ is invariant closed. 
Since $F = X - U$ and $\hat{\F}(x) \subseteq U$, we have $F \cap \hat{\F}(x) = \emptyset$ and so $\hat{\F}(F) \cap \hat{\F}(x) = \emptyset$. 
Therefore the complement $X - \hat{\F}(F)$ is an invariant open \nbd of $\hat{\F}(x)$ contained in $U$. 
This means that the class decomposition $\hat{\F}$ is weakly upper semi-continuous. 

Conversely, suppose that $\F$ is pointwise almost periodic and the class decomposition $\hat{\F}$ is weakly upper semi-continuous. 
Fix a closed subset $F \subseteq X$. 
We claim that the complement $X - \hat{\F}$ is open. 
Indeed, fix a point $x \in X - \hat{\F}$. 
Then the complement $U := X - F$ is an open \nbd of $\hat{\F}(x)$. 
By weak upper semi-continuity of $\hat{\F}$, there is an open $\hat{\F}$-invariant \nbd $V \subseteq U$ of $\hat{\F}(x)$. 
Since $V \subseteq U = X - F$, we have $V \cap F = \emptyset$. 
The $\hat{\F}$-invariance of $V$ implies that $V \cap \hat{\F}(F) = \emptyset$ and so $V \subseteq X -  - \hat{\F}$. 
This means that the complement $X - \hat{\F}$ is open. 
Therefore $\F$ is weakly almost periodic. 
\end{proof}

Moreover, we have the following implication. 

\begin{lemma}\label{lem:wap_r02}
The class space of a weakly almost periodic semi-decomposition $\F$ on a normal space $X$ is Hausdorff. 
Moreover, if $X$ is Hausdorff, then $\F$ is $R$-closed. 
\end{lemma}

\begin{proof}
Let $\F$ be a weakly almost periodic semi-decomposition on a normal space $X$. 
Since $\F$ is pointwise almost periodic, 
the class semi-decomposition is a decomposition $\hat{\F} = \{\overline{\F(x)} \mid x \in X\}$ and so the element closure relation of $\hat{\F}$ is an equivalence relation on $X$. 
Since the element closure relation $R$ of $\F$ is also one of the decomposition $\hat{\F}$, the element closure relation of $\F$ is the same equivalence relation on $X$. 
Fix two distinct classes $L \neq L' \in \hat{\F}$. 
By normality of $X$, there are disjoint open \nbds $U, U'$ of $L$ and $L'$ respectively. 
Since the complements $F := X - U$ and $F' := X - U'$ are closed, the saturation $\hat{\F}(F)$ and  $\hat{\F}(F')$ are closed. 
By construction, we have $\hat{\F}(F) \cap L = \emptyset$ and $\hat{\F}(F') \cap L' = \emptyset$. 
Then $L \subseteq X - \hat{\F}(F) \subseteq U$ and $ L' \subseteq X - \hat{\F}(F') \subseteq U'$. 
This means that $X - \hat{\F}(F)$ and $X - \hat{\F}(F')$ are disjoint open $\hat{\F}$-invariant \nbds of $L$ and $L'$ respectively. 

Suppose that $X$ is Hausdorff. 
Since $\hat{\F} = \{\overline{\F(x)} \mid x \in X\}$, the element closure relation $R$ of $\F$ is also one of the decomposition $\hat{\F}$. 
\cite[Proposition~8.3.8 p.79]{bourbaki2013general} implies that the element closure relation $R$ of $\F$ is closed. 
\end{proof}

\subsubsection{Properties of minimality for semi-decompositions}

We have the following statement. 

\begin{proposition}\label{th:mini}
Any minimal semi-decomposition on a topological space is of characteristic $0$, $R$-closed, and weakly almost periodic. 
\end{proposition}

\begin{proof}
Let $\F$ be a minimal semi-decomposition on a topological space $X$. 
Fix a point $x \in X$. 
By Lemma~\ref{lem:prolong}, we have $X = \overline{\F(x)} \subseteq D(x)$ and so $\overline{\F(x)} = D(x)$. 
This means that $\F$ is of characteristic $0$. 
By definition of the element closure relation $R$, the minimality implies that $R = \{ (x,y) \in X \times X \mid y \in \overline{\F(x)} \}= X \times X$ is closed. 
Therefore $\F$ is $R$-closed. 
From the minimality, the class semi-decomposition $\hat{\F}$ consists of one element $X$ and so the semi-decomposition $\F$ is pointwise almost periodic. 
Since the union $\bigcup_{x\in A} \overline{\F(x)}$ of element closures for any non-empty closed subset $A \subseteq X$ is $X$ and so closed. 
This means that $\F$ is weakly almost periodic. 
\end{proof}

\subsection{Relations of variants of recurrence for decompositions}

We have the following statement. 

\begin{lemma}\label{lem:wap}
Any $R$-closed decomposition with compact element closure on a Hausdorff space is weakly almost periodic.
\end{lemma}

\begin{proof}
Let $\F$ be an $R$-closed decomposition with compact element closure on a Hausdorff space $X$. 
The proof of \cite[Lemma~1.5]{yokoyama2013recurrence} is true when the Hausdorff separation axiom for $X$ is replaced by the compact element closure condition. 
Therefore the assertion hold. 
\end{proof}

We have the following equivalence.
 
\begin{lemma}\label{lem:02}
The following statements are equivalent for a decomposition on a $T_3$ space:
\\
{\rm(1)} The decomposition is of characteristic $0$.  
\\
{\rm(2)} The decomposition is $R$-closed.
\\
{\rm(3)} The decomposition is pointwise almost periodic and the class space is Hausdorff. 
\\
Moreover, if any element closures are compact, then any equivalent conditions {\rm(1)}--{\rm(3)} imply the weakly almost periodicity of the group action. 
\end{lemma}

\begin{proof}
By \cite[Corollary 1.4 and Lemma 2.2]{yokoyama2013recurrence}, the assertions {\rm(2)} and {\rm(3)} are equivalent. 
From \cite[Corollary 3.1]{yokoyama2013recurrence}, the assertions {\rm(1)} and {\rm(2)} are equivalent. 
Consider a decomposition on a compact Hausdorff space $X$. 
By Lemma~\ref{lem:wap}, if any element closures are compact, then the group action is weakly almost periodic. 
\end{proof}



We have the following proof. 

\begin{proof}[Proof of Theorem~\ref{th:equi}]
Since any metrizable space is $T_3$, Corollary~\ref{lem:equi001} and Lemma~\ref{lem:02} imply the assertions {\rm(1)}--{\rm(4)}. 
By Lemma~\ref{lem:wap_r02}, the assertion {\rm(4)} implies the assertion {\rm(5)}. 
From \cite[Lemma 3.6]{yokoyama2017low}, the assertion {\rm(4)} implies the assertion {\rm(6)}. 
\end{proof}

Theorem~\ref{th:equi} implies the following statement.

\begin{corollary}\label{lem:equi2R_metric}
Every equicontinuous group action on a metrizable space is of characteristic $0$ and is $R$-closed, and the orbit class space is Hausdorff.  
Moreover, if any orbit closures are compact, then the action is weakly almost periodic. 
\end{corollary}

\section{Incomparability among variants of recurrence for actions of groups and semi-groups}

We recall some concepts to show incomparability. 

\subsection{Concepts for {\rm(}semi-{\rm)}group actions}

\subsubsection{{\rm(}Semi-{\rm)}group actions}
A continuous mapping $v \colon T \times X \to X$ on a topological space $X$ and a topological semi-group $T$ with unit element is a {\bf semi-group action} of $T$ (or a $T$-action) on $X$ if the restriction $v(t,\cdot) \colon X \to X$ for any $t \in T$ is homeomorphic such that the restriction $v(e,\cdot)$ for the unit element $e \in T$ is an identity mapping on $X$ and $v(g,v(h,x)) = v(gh,x)$ for any $g,h \in T$.
A semi-group action is a {\bf group action} if the semi-group is a group. 
%
The subset $\{g(x) \mid g \in T \}$ for any point $x \in X$ is called the {\bf orbit} of $x$ and is denoted by $\bm{T(x)}$. 
The closure of the orbit of a point is called the orbit closure of the point. 
For a subset $A \subseteq X$, put $T(A) := \bigcup_{x \in A}T(x)$. 
A subset $A \subseteq X$ is ($T$-)invariant if $A = T(A)$. 

By a {\bf flow}, we mean a continuous $\R$-action. 
A point $x$ of a flow $v \colon \R \times X \to X$ is {\bf recurrent} if $x \in \omega(x) \cup \alpha(x)$, where $\omega(x) := \bigcap_{n\in \mathbb{R}}\overline{\{v(t,x) \mid t > n\}}$ is the $\omega$-limit set and $\alpha(x) := \bigcap_{n\in \mathbb{R}}\overline{\{v(t,x) \mid t < n\}}$ is the $\alpha$-limit set. 
A flow is {\bf recurrent} if each point is recurrent. 
Notice that a $\Z$-action can be identified with a homeomorphism canonically. 

From now on, we use a $T$-action or a semi-group action of $T$ for a semi-group action of a topological semigroup $T$ with unit element. 
For any $T$-action on a topological space $X$, consider the induced $T$-action on $X \times X$ defined by $g(x,y) := (g(x), g(y))$. 
More generally, the product action $G^{\kappa}$ on the product space $X^\kappa$ is defined by $g(x_\alpha) := (g(x_\alpha)) \in X^\kappa$. 



\subsubsection{Pointwise almost periodicity, $R$-closedness, characteristic $0$, weakly almost periodicity for group actions}

Denote by $\F_T$ the set of orbits of a semi-group action of $T$. 

\begin{definition}
A semi-group action of $T$ is {\bf pointwise almost periodic} {\rm(resp.} {\bf $R$-closed}, {\bf characteristic $0$}, {\bf weakly almost periodic}{\rm)} if the semi-decomposition  $\F_T$ is pointwise almost periodic {\rm(resp.} $R$-closed, characteristic $0$, weakly almost periodic{\rm)}. 
\end{definition}

Notice that a semi-group action on a topological space $X$ is pointwise almost periodic if and only if the orbit class space is $T_1$. 
Though this definition of pointwise almost periodicity for a group action is different from one in \cite[p.35]{auslander1988minimal}, these definitions coincide when the semigroup on the topological space is a group on a locally compact Hausdorff space \cite[Theorem~7 p.11 and Lemma~3 p.37]{auslander1988minimal}. 
When the space is compact, Gottschalk's original definition of weakly almost periodicity for group actions corresponds to our definition introduced by Mason (see \cite[Theorem~4.24 p.34]{gottschalk1955topological} for details). 
There are pointwise almost periodic group actions that are not of characteristic $0$\cite{auslander2007recurrence}. 

Moreover, there is a pointwise recurrent smooth flow on a torus that is not pointwise almost periodic. 
Indeed, consider an irrational rotating vector field $X_0(x,y) = (1, \theta)$ on a torus $\mathbb{T}^2 = (\R/\Z)^2$ for any irrational number $\theta$. 
Fix a point $p_0 \in \mathbb{T}^2$. 
Using a smooth bump function $f \colon \mathbb{T}^2 \to [0,1]$ with $f^{-1}(1) = \{ p_0 \}$, the resulting vector field $(1 -f) X_0$ generates a flow $v$ on $\mathbb{T}^2$ with $\mathop{\mathrm{Sing}}(v) = \{ p_0 \}$ such that any orbits except $\{ p_0 \}$ are dense.  
This means that the flow $v$ is a smooth recurrent flow but is not pointwise almost periodic. 

In the low-dimensional cases, the characteristic $0$ property coincides with other concepts. 
For instance, in \cite{auslander2007recurrence,hattab2011pointwise}, it is shown that the following are equivalent for each action of a finitely generated group $G$ on either a compact zero-dimensional space or a graph $X$: (1) Pointwise recurrent;  (2) Pointwise almost periodic; (3) Characteristic $0$. 
Moreover, the equivalence holds for equicontinuously generated groups on locally compact zero-dimensional spaces \cite{reid2020equicontinuity}. 
In addition, the following are equivalent for a countable group $G$ on a local dendrite:  (1) Pointwise almost periodic;  (2) Hausdorff for the orbit class space; (3) Characteristic $0$ \cite{Marzougui2019min}.

\subsubsection{Orbit classes and orbit class spaces for semi-group actions}
Let $v: T \times X \to X$ a semi-group action of $T$. 
An element class containing $x$ of the semi-decomposition $\F_T$ is called the {\bf orbit class} of $x$.  
If the class semi-decomposition $\hat{\F}_T$ of the semi-decomposition $\F_T$ is a decomposition, then the decomposition space of $\hat{\F}_T$ is called the {\bf orbit class space} and denoted by $\bm{X/\hat{T}}$. 

\subsubsection{Equicontinuity}
We equip the topology of the uniformity with any uniform space. 

\begin{definition}
A semi-group action of $T$ on a uniform space $(X,\mathcal{U})$ is {\bf equicontinuous} if, for any $\alpha \in \mathcal{U}$ there is an entourage $\beta \in \mathcal{U}$ with $T (\beta) \subseteq \alpha$ {\rm(i.e.} $(gx, gy) \in \alpha$ for any $(x,y) \in \beta$ and any $g \in T${\rm)}. 
\end{definition}

A semi-group action on a uniformizable space $(X,\mathcal{O})$ is {\bf equicontinuous} if there is a uniformity $\mathcal{U}$ on $X$ that induces the topology $\mathcal{O}$ and that the semi-group action on $(X,\mathcal{U})$ is equicontinuous. 
When $X$ has a compatible metric $d$ to the topology of the uniformity $\mathcal{U}$, a $T$-action on $X$ is equicontinuous if and only if, for any $\varepsilon \in \R_{>0}$ there is a positive number $\delta \in \R_{>0}$ such that $d(g(x), g(y)) < \varepsilon$ for any points $x, y \in X$ with $d(x,y) < \delta$ and $g \in T$. 
In other words, a $T$-action on a metric space $(X,d)$ is equicontinuous if and only if, for any $\varepsilon \in \R_{>0}$ there is a positive number $\delta \in \R_{>0}$ with $\sup \{ d(g(x), g(y))\mid g \in T, x, y \in X, d(x,y) < \delta \} < \varepsilon$. 
%
A semi-group action on a metrizable space $X$ is {\bf equicontinuous} if there is a compatible metric $d \colon X \times X \to \R$ such that the semi-group action on the metric space $(X,d)$ is equicontinuous. 
By definition of Hausdorff distance, equicontinuity of semi-group actions implies weakly equicontinuity for the semi-decomposition of orbits. 
Moreover, the suspension construction for group actions (cf. \cite[Chapter~3]{CC2003I}) in some settings induces the equicontinuity of the induced foliations from one of group actions (see details for \cite[2B]{dyer2017molino}).

%
Furthermore, the equivalence of equicontinuity and the closedness of the orbit closure relation does not hold \cite{knight1986cer,salem2018group}, even for distal homeomorphisms (i.e. continuous $\Z$-action) on a one-dimensional compact metric space (see \cite[Introduction]{salem2018group} and an example in \S~\ref{sec:ex_r_nonec} for details). 
In addition, there are non-equicontinuous $R$-closed discrete and continuous dynamical systems \cite{elaydi1982char,elaydi1990variation}.
Therefore, the proofs of \cite[Proposition 4.8 and Theorem~1.2]{Hattab2017flow}, and \cite[Theorem 2.1 and Propositions 4.2--4.3]{salem2017dendrite} do not work. 

\subsubsection{Distal property}

Distal property, introduced by Hilbert \cite{Hilbert1956geomtry} to topologically characterize a rigid group of motions, is defined as follows. 
\begin{definition}
A semi-group action of $T$ on a topological space $X$ is {\bf distal} {\rm(}cf. \cite[Definition~4.5 p.41]{ellis2014automorphisms} {\rm)} if $P(X) = \Delta$, where $\Delta := \{ (x,x) \in X \times X \mid x \in X \}$ and $P(X) := \{ (x,y) \in X \times X \mid \overline{T(x,y)} \cap \Delta \neq \emptyset \}$. 
\end{definition}
%
When $X$ has a compatible metric $d$ to the topology, a $T$-action on $X$ is distal if and only if $x = y$ for any triple $x,y,z \in X$ and any sequence $(g_n)_{n \in \Z_{>0}}$ with $z = \lim_{n \to \infty} g_n(x) = \lim_{n \to \infty} g_n(y)$. 

Notice that the definition of distal property for group action on topology spaces coincides with one for uniform spaces under the compact Hausdorff condition. 
%
Note that there is a non-distal minimal group action (see \cite[an example after Theorem~2.9 in pp.392--393]{keynes1967study}. 
Though every equicontinuous surjective semiflow (i.e. action by a semi-group with identity) on a compact Hausdorff space is distal \cite[Theorem~2.1]{Auslander2019mini}, notice that this implication for semiflows on compact Hausdorff spaces does not hold (see \cite[Example~1.2]{Auslander2019mini}). 
In addition, there is a countable non-equicontinuous minimal $R$-closed distal group action on a zero-dimensional compact metric space $X$ in \cite[Example 1--2. \S~3.4 pp.265--266]{mcmahon1976connectedness} as mention in \cite[p.108]{auslander2007recurrence} (cf. \cite[Proposition 4.1.]{liang2021non}). 

\subsubsection{Regularity}

Consider a semi-group action of $T$ on a topological space. 
A point $x \in X$ is {\bf regular} \cite{Kaul1973irr} if for any open subset $U$ and any subset $P \subseteq T$ with $\overline{P(x)} \subseteq U$, there is an open \nbd $V$ of $x$ such that $\bigcup_{y \in V} P(y) \subseteq U$. 
A semi-decomposition on a topological space $X$ is of {\bf regular} if so is any point of $X$. 
The regularity and equicontinuity coincide for group actions with compact orbit closure on locally compact Hausdorff spaces \cite{Kaul1977zero} (cf. \cite[Theorem 2.9]{elaydi1981weakly}). 

\subsection{Non-equicontinuous group actions}\label{sec:ex_nonec}

We construct examples to show that $R$-closedness does not imply equicontinuity.

\subsubsection{A non-equicontinuous distal $R$-closed weakly almost periodic pointwise periodic analytic flow on a closed surface}\label{sec:ex_analytic}

By an analogous construction to \cite[Example~3.1]{elaydi1990variation}, we construct a non-equicontinuous distal $R$-closed weakly almost periodic pointwise periodic analytic flow on a torus $\mathbb{T}^2 = (\R/\Z)^2$ as follows. 
Define a toral flow $v \colon \R \times \mathbb{T}^2 \to \mathbb{T}^2$ by $v(t, [x,y]) := [x, y + (2+ \cos x) t]$. 
Then the orbit class space is a circle, and so the action is $R$-closed. 
Moreover, any orbits are periodic and so are homeomorphic to circles. 
This implies distal property. 
For any point $p_1 = [x_1, y_1] \in \mathbb{T}^2$ and for any positive number $\delta$, there is a positive number $t \in \R$ such that $d(v(t,p_1), v(t,p_2)) > 1/3$, where $p_2 := [x_1 + \delta/2, y_2]$. 
This means that $v$ is not equicontinuous. 

Furthermore, the time-one map is a non-equicontinuous non-$R$-closed distal pointwise periodic analytic diffeomorphism on a torus because the orbit class space is not Hausdorff.

%

\subsubsection{Non-weak equicontinuity, distal property, and $R$-closedness to any countable product actions for group actions on non-compact paracompact manifolds}\label{sec:ex_r_nonec}
%
There is a non-weakly equicontinuous distal action on a zero-dimensional paracompact manifolds $\mathcal{C}$ such that the orbit closure relation on any countable product $\mathcal{C}^{\kappa}$ is closed. 
%
Indeed, consider a metric $d \colon \R \to \R$ defined by $d(x,y) = |x^3 - y^3|$ on $\R$ and 
a diffeomorphism $f \colon \R \to \R$ by $f(x) = x+2$. 
By construction, the group $\Gamma$ generated by one generator $f$ is not weakly equicontinuous. 
Moreover, any orbits of $\Gamma$ are non-compact closed subsets.  
Equip the countable product $\R^{\kappa}$ with a product metric (e.g. $d_{\kappa}((x_i)_i,(x'_i)_i) = \sum_{i \geq 1} d(x_i, x'_i)/2^i$). 
The product $\Gamma$-action $v^{\kappa} \colon \Gamma \times \R^{\kappa} \to \R^{\kappa}$ defined by $v^{\kappa}(f^n, (x_i)_i) = (x_i + 2n)_i \in \R^{\kappa}$ satisfies that the orbit closure relation on $\R^{\kappa}$ is closed.
Denote by $\mathcal{C}_0$ the Cantor ternary set on $[0,1]$. 
Put $\mathcal{C} := \bigsqcup_{n \in \Z} f^n(\mathcal{C}_0)$. 
Put $g := f|_{\mathcal{C}}$. 
By construction, the group $G$ generated by one generator $g$ is non-weakly equicontinuous. 
Since the closedness of the orbit closure relation is invariant under taking restrictions of domains, the  group action $v^{\kappa}| \colon G \times \mathcal{C}^{\kappa} \to \mathcal{C}^{\kappa}$ satisfies that the orbit closure relation on $\mathcal{C}^{\kappa}$ is closed for any cardinal $\kappa$.

\subsection{Incomparability among variants of recurrence for group actions}\label{sec:ex}

Theorem~\ref{th:incomparable_semi}(1) follows from the examples below in \S\ref{sec:ex02}--\ref{sec:ex03}. 
Theorem~\ref{th:incomparable_semi}(5) follows from the examples below in \S\ref{sec:ex03}--\ref{sec:ex_ec_nonmini}.

%
%
%
%


\subsubsection{Non-$R$-closed distal smooth flow on a sphere}\label{sec:ex02}

There is a non-$R$-closed (non-weakly equicontinuous) distal (pointwise almost periodic) smooth flow on a sphere $\mathbb{S}^2$. 
Consider a spherical rotating vector field with exactly two critical points $x_N, x_S$. 
Using a bump function, replace the periodic orbit $\gamma$ on the equator into singular points. 
Then the flow $v$ generated by the resulting vector field consists of closed orbits such that $\mathop{\mathrm{Sing}}(v) = \{ x_N, x_S \} \sqcup \gamma$, where $\mathop{\mathrm{Sing}}(v)$ is the set of singular points of $v$. 
Notice that $\mathop{\mathrm{Per}}(v) = \mathbb{S}^2 - \mathop{\mathrm{Sing}}(v)$. 
Moreover, the restriction of $v$ to the periodic orbit, which is homeomorphic to a circle, is a rotation on the circle up to topological equivalence. 
This means that the restriction of $v$ to any orbit is either a trivial flow on a critical point or a rotation on a circle up to topological equivalence. 
Therefore the flow $v$ is distal. 
Since the orbit class space $\mathbb{S}^2/\hat{v}$ can be obtained by collapsing periodic orbits into singletons, it is not Hausdorff and so $v$ is not $R$-closed. 
By Theorem~\ref{th:equi}, the flow $v$ is neither equicontinuous nor weakly equicontinuous. 

\subsubsection{$R$-closed weakly almost periodic non-distal non-weakly equicontinuous minimal action on a Cantor set}\label{sec:ex03}

The restriction to the Cantor minimal set of the Denjoy diffeomorphism is a non-distal non-weakly equicontinuous weakly almost periodic $R$-closed minimal homeomorphism on a Cantor set. 
Indeed, recall the Denjoy construction. 
Fix an irrational number $r \in \R - \Q$, a translation $\widetilde{f_0} \colon \R \to \R$ by $\widetilde{f_0} (x) := x +r$, and an irrational rotation $f_0 \colon \mathbb{S}^1 \to \mathbb{S}^1$ on a circle $\mathbb{S}^1 := \R/\Z$ by $f_0([x]) := [x+r]$. 
Since $f_0$ is minimal, there is an increasing sequence $(n_i)_{i \in \Z_{>0}}$ such that the sequence $(f_0^{n_i}(0))_{i \in \Z_{>0}}$ converges to $[1/2] \in \mathbb{S}^1$. 
Consider a sequence $\{l_n := 1/(1+n^2) \mid n \in \Z \}$ and a sequence $\{I_n \mid n \in \Z \}$ of half-open intervals $I_n$ with the length $l_n$.  
Put $L := \sum_{n \in \Z} l_n = \sum_{n \in \Z} 1/(1+n^2) < \infty$ and define $y_n := \widetilde{f_0}^n(0) - \lfloor \widetilde{f_0}^n(0) \rfloor \in [0,1)$. 
Insert the interval $I_n$ at the point $y_n \in [0,1)$ and denote by $S^1 := \R/(1+L)\Z$ the resulting circle and by $f \colon S^1 \to S^1$ the resulting homeomorphism with a unique minimal set $\mathcal{C}$, which is a Cantor set. 
By construction, we have that $\mathcal{C} = S^1 - \bigcup_{n \in \Z} I_n$ and that the restriction $g:= f|_{\mathcal{C}}$ is a minimal homeomorphism on a Cantor set. 
Therefore $g$ is $R$-closed and weakly almost periodic. 
Collapsing the intervals $I_n$ into singletons, we obtain a semi-conjugate $h \colon S^1 \to \mathbb{S}^1$ from $f$ to $f_0$ (i.e. $f \circ h = h \circ f_0$). 
By construction, the inverse image $h^{-1}([1/2]) \subset \mathcal{C}$ is a singleton and is denoted by $p$. 
Define points $0_-, 0_+$ by $\{ 0_-, 0_+ \} = \partial I_0 \subset \mathbb{S}^1$. 
By $\lim_{n \to \infty} l_n = \lim_{n \to \infty} 1/(1+n^2) = 0$ and $\lim_{n \to \infty} h(f^{n_i}(0_-)) = \lim_{n \to \infty} h(f^{n_i}(0_+)) = [1/2]$, since the inverse image $h^{-1}([1/2])$ is a singleton, we obtain $\lim_{n \to \infty} f^{n_i}(0_-) = \lim_{n \to \infty} f^{n_i}(0_+) = h^{-1}([1/2]) \in \mathcal {C}$. 
This means that $g = f|_{\mathcal{C}}$ is neither distal nor weakly equicontinuous. 

\subsubsection{Equicontinuous non-minimal actions}\label{sec:ex_ec_nonmini}
Any trivial actions on metric spaces which contains at least two points is not minimal but equicontinuous. 
In particular, the actions are weakly equicontinuous and distal.

\subsection{Incomparability among variants of recurrence for semi-group actions}

Theorem~\ref{th:incomparable_semi}(2) follows from the example in \S\ref{sec:ex02} and one below in \S\ref{sec:ec_npap}. 
Incomparability between distal property and pointwise almost periodicity follows from the example in \S\ref{sec:ex03} and one below in \S\ref{sec:ec_npap}. 
Other incomparability in Theorem~\ref{th:incomparable_semi}(3) follows from the example in \S\ref{sec:ex02} and one below in \S\ref{ex:r_non_pap}. 
Theorem~\ref{th:incomparable_semi}(4) follows from the example in \S\ref{sec:ex03} and one below in \S\ref{ex:r_non_pap}. 
%
%
%
%
%
%

\subsubsection{Equicontinuous $R$-closed distal semi-group action which is not pointwise almost periodic}\label{ex:r_non_pap}
There is an $R$-closed equicontinuous distal semi-group action which is not pointwise almost periodic such that the orbit closure relation is not pointwise almost periodic. 
Indeed, let $X = \{ 0, 1 \}$ be the discrete metric space and $f \colon X \to X$ a mapping defined by $f(0) = f(1) = 0$. 
Then the orbit closure relation of the semi-group generated by $f$ is $R = \{ (x,y) \in X \times X \mid y \in \overline{\F(x)} \} = \{(0,0), (1,0), (1,1) \}$ is not symmetric. 
Since the orbit closure of $1$ is not minimal, the action is not pointwise almost periodic. 
Moreover, we have $D(0) = \{ 0 \} \neq X = D(1)$.  
In addition, since the singleton of any point is open, the action is regular.

\subsubsection{Equicontinuity non-distal non-pointwise almost periodicity for contractions}\label{sec:ec_npap}
Recall that a self-mapping $f \colon X \to X$ on a (possibly compact) metric space $(X,d)$ is a contraction if there is a constant $C \in (0,1)$ such that $(d(f(x), f(y)) \leq C d(x, y)$ for any points $x, y \in X$. 
%
Then any semi-group action generated by a non-identical contraction with a unique fixed point is equicontinuous but neither distal, pointwise almost periodic, nor weakly almost periodic. 
%
Indeed, by definition of equicontinuity and distal, the semi-action is equicontinuous but not distal. 
For any orbit $O$ and any points $x \neq y \in O$, the contraction implies that $\overline{O(x)} \neq \overline{O(y)}$. 
This means that any orbits which are not singletons are not minimal and so that the set of orbit closures is not a decomposition. 
Therefore the action is neither pointwise almost periodic nor weakly almost periodic. 
By Corollary~\ref{lem:equi001}, the equicontinuous action is $R$-closed and of characteristic $0$.  


\vspace{10pt}

{\bf Acknowledgement}: 
The author wishes to thank Hiraku Nozawa for posing a motivative question.

\bibliographystyle{abbrv}
\bibliography{yt20211011}

\end{document}